\documentclass[a4paper,12pt,onecolumn]{IEEEtran}

\IEEEoverridecommandlockouts    

\usepackage{mathptmx} 
\usepackage{times} 
\usepackage{amsmath} 
\usepackage{amssymb}  
\usepackage{amsxtra}
\usepackage{amsfonts}

\usepackage{mathrsfs}

\usepackage{graphicx}

\usepackage{color}
\usepackage{tikz}

\DeclareSymbolFont{bbold}{U}{bbold}{m}{n}
\DeclareSymbolFontAlphabet{\mathbbm}{bbold}

\DeclareMathAlphabet{\mathcal}{OMS}{cmsy}{m}{n}

\newcommand{\norm}[1]{\ensuremath{\left\| #1 \right\|}}

\newcommand{\Real}{\ensuremath{\mathbb R}}

\newcommand{\LT}{\mathbf{L}_{\rm 2}}
\newcommand{\LTE}{\mathbf{L}_{\rm 2e}}
\newcommand{\jjmath}{\mathfrak{j}}
\newcommand{\jw}{\jjmath\omega}

\newcommand{\SSS}{\mathscr{S}}
\newcommand{\HHH}{\mathscr{H}}

\newcommand{\RRR}{\mathscr{R}}

\newcommand*{\dif}{\mathop{}\!\mathrm{d}}

\newtheorem{theorem}{Theorem}

\newtheorem{proposition}{Proposition}
\newtheorem{lemma}{Lemma}
\newtheorem{remark}{Remark}

\title{\LARGE \bf  Integral quadratic constraints for asynchronous
  sample-and-hold links}

\author{Michael Cantoni,
  Chung-Yao Kao, and
  Mark A.~Fabbro
  \thanks{*Supported in part by the Australian Research Council (DP130104510) and Ministry of Science and Technology of Taiwan (MOST
    106-2221-E-110-007-MY2).}
  \thanks{M.~Cantoni and M.A.~Fabbro are with the Department of Electrical and Electronic
    Engineering, The University of Melbourne, Parville, VIC 3010,
    Australia. {\tt\footnotesize cantoni@unimelb.edu.au,
      markfabbro7@gmail.com}}%
  \thanks{C.-Y.~Kao is with the Department of Electrical Engineering, National
    Sun Yat-Sen University, Kaohsiung, 80424,
    Taiwan. {\tt\footnotesize cykao@mail.nsysu.edu.tw} }%
}

\begin{document}

\maketitle

\thispagestyle{empty}
\pagestyle{empty}

\begin{abstract}
  A model is proposed for a class of asynchronous sample-and-hold operators that is relevant in the analysis of embedded and networked systems. The model is parametrized by characteristics of the corresponding time-varying input-output delay. Uncertainty in the relationship between the timing of zero-order-hold update events at the output and the possibly aperiodic sampling events at the input means that the delay does not always reset to a fixed value. This is distinct from the well-studied synchronous case in which the delay intermittently resets to zero at output update times. The main result provides a family of integral quadratic constraints that covers the proposed model. To demonstrate an application of this result, robust $\LT$ stability and performance certificates are devised for an asynchronous sampled-data implementation of a feedback loop around given linear time-invariant continuous-time open-loop dynamics. Numerical examples are also presented.
\end{abstract}
\begin{IEEEkeywords}
Integral quadratic constraints, robustness analysis, sampled-data networks, time-varying delay
\end{IEEEkeywords}

\section{Introduction}

The digital implementations of controllers generically involve plant output sampling and plant input updates at discrete time instants. This leads to time-varying closed-loop dynamics, even when the plant and controller are time invariant. In particular, at times between updates, the plant input is held constant at a value determined by a plant output sample taken at a varying time in the past~\cite{frid04}. 

A well studied approach to the analysis of sampled-data control
systems, called the ``input delay'' approach, applies when the
discrete-time controller implementation commutes with the hold
operation, as in the case of static gain feedback. Combinding the
sampling and the hold operations
yields a closed-loop that can be modelled in a
purely continuous-time fashion. More generally, sample-and-hold
operators arise in the study of digital networks of continuous-time
dynamical systems~\cite{FabbroThesis}. It is well known that
sample-and-hold operators with synchronized
sampling and hold event timing (i.e., synchronous sample-and-hold) can
be modelled as a saw-tooth time-varying delay that resets to zero at
the possibly non-uniform update/sampling
instants~\cite{frid04},~\cite{Mir07}--\cite{Kao16}.

Implementation resource limitations may lead to asynchronous hold update and sample event sequences at the output and input, respectively, of a sample-and-hold link. This can arise
when there is variability 
in the
time consumed by the mechanisms used to process samples from sensors and communicate data to actuators. 
When the timing of sample events and hold update events is not synchronized
(i.e., asynchronous sample-and-hold),
the result is a time-varying delay that does not always reset to zero (or another constant value) at hold update events. But like the synchronous case with aperiodic sampling, exact description of the time-varying delay may not be possible ahead of time. Instead, the delay can be considered uncertain and abstractly modelled by bounds on inter-update intervals, inter-sample intervals, and other relationships reflecting asynchrony between output and input events. Ultimately, the model developed below is expressed in terms of a perturbation of the identity (i.e., relative to the ideal link).

The main contribution in what follows is a set of results that establish gain bound and passivity properties of the aforementioned perturbation of the identity. This extends results in~ and~\cite{Mir07}, \cite{Fuj09} and~\cite{Kao16}, where the synchronous special case is studied; also see~\cite{KaoLincoln}--\cite{Rahn18}. Combining the gain and passivity  properties yields a family of integral quadratic constraints (IQCs) that cover the uncertain asynchronous sample-and-hold operator in a way that is suitable for robustness analysis in the vein of~\cite{MegRan}. The results established here go beyond the preliminary versions of the work reported in~\cite{CanFabKao18} and \cite{CanFabKao19}. In particular, an additional characteristic of the time-varying delay is accounted for more explicitly here. Its consideration leads to a less conservative IQC cover, as discussed further in Sections~\ref{sec:IQCcovers}. To illustrate use of the proposed model and its IQC characterization, robust stability and performance certificates are devised for the sampled-data implementation of a feedback loop, subject to resource limitations that lead to asynchronous sampling and update events. 

The input-output context of the IQC approach followed here is a distinguishing feature relative to the more common approaches to sampled-data system analysis based on hybrid/impulsive state-space modelling~\cite{Hes07}--\cite{Het17}. Further, in the state-space literature it is standard to relate all sample and update events to a single time sequence, for which one interval bound holds. By contrast, in the structured input-output approach developed here, characteristics of the sample and update event sequences are bounded individually, and with respect to 
each other. Recent related work on asyncrohonous sample-data systems
within a state-space context can be found
in~\cite{Fiacc16}--\cite{Etienne19}. While gain bounds on related
operators play a role in some of these papers, the combined
exploitation of gain and passivity properties is not considered. It is
shown by example here that the consideration of both together can
reduce conservativeness in the resulting analysis.


The rest of the paper is organized as follows. First some notation and terminology are established. In Section~\ref{sec:asynchSH}, the class of uncertain asynchronous sample-and-hold operators is formally defined. The aforementioned gain bound and passivity properties are then derived, and combined to formulate the IQC covers of the asynchronous sample-and-hol model, in Section~\ref{sec:IQCcovers}. Corresponding robust stability and performance certificates for the example sampled-data feedback loop are devised in Section~\ref{sec:examples}, where numerical result are also presented. Concluding remarks are provided in Section~\ref{sec:conc}.

\section{Preliminaries}

The non-negative integers and the reals are denoted by $\mathbb{N}_0$ and $\mathbb{R}$, respectively, and
$\mathbb{R}_{>0}=(0,+\infty)\subset\mathbb{R}$. The space of square integrable functions defined on $\mathbb{R}_{\geq 0} =[0,+\infty)\subset\mathbb{R}$ is denoted by $\LT$, and the usual $\LT$ norm and inner product are denoted by $\|\cdot\|_{\LT}$ and $\langle\cdot,\cdot\rangle_{\LT}$, respectively.
The extended $\LT$ space is denoted by $\LTE$. This consists of
functions $f:\mathbb{R}_{\geq 0}\rightarrow\mathbb{R}$ that satisfy $P_\tau f\in\LT$ for
$\tau>0$, where $P_\tau$ is the truncation
operator; i.e., $(P_\tau f)(t)=f(t)$ for $t\leq \tau$, otherwise
$(P_\tau f)(t)=0$. 
The set of right continuous functions defined on $\mathbb{R}_{\geq 0}$ is denoted by $\mathbf{C}_r$. 


Let $G:\LTE\rightarrow\LTE$ be given. This is called a bounded operator if $u\in\LT \implies Gu\in\LT$ and $\|G\|=\sup_{u\in\LT}\|Gu\|_{\LT}/\|u\|_{\LT}$ is finite. If
$P_\tau G P_\tau - P_\tau G
= 0$
for all $\tau > 0$, then $G$ is called causal, and when it is also bounded, $G$ is called stable. If
$G$ is linear and bounded (not necessarily causal), the adjoint of the restriction to $\LT$ is denoted by $G^*$; i.e., $\langle v, Gu \rangle_{\LT}=\langle G^*v,u\rangle_{\LT}$ for $u,v\in\LT$.
%
If $G$ is stable, linear, and time-invariant (i.e., commutes with forward shift), then its restriction to $\LT$ corresponds to multiplication by a frequency domain transfer function. This transfer function, also denoted $G$ for convenience, is analytic in the right-half plane, with
$\|G\|=
\mathrm{ess}\sup_{\omega \in \Real}|G(\jw)|$,
where $\jjmath=\sqrt{-1}$.
If $G$
admits the rational transfer function $G(s)=C(sI-A)^{-1}B+D$, the collection of matrices $(A,B,C,D)$ is called a state-space realization. The identity matrix is denoted by $I$.

Let $\Delta:\LTE\to\LTE$ be stable. Given self-adjoint $\Pi:\LT\rightarrow\LT$, the bounded causal operator
$\Delta:\LT\to\LT$ is said to satisfy the IQC defined by $\Pi$ if 
\begin{align*}
    \left\langle
    \begin{bmatrix} v \\ \Delta v
    \end{bmatrix}, 
    \Pi 
    \begin{bmatrix} v \\ \Delta v \end{bmatrix}
    \right\rangle_{\LT} \geq 0 ~\text{ for all } v\in\LT. 
\end{align*}
When this property holds, it is written that $\Delta\in\mathrm{IQC}(\Pi)$.
%
Dependence of the so-called multiplier $\Pi$ on a parameter $X$ is denoted by $\Pi(X)$.

Given causal $G:\LTE \rightarrow\LTE$ and
$\Delta:\LTE \rightarrow\LTE$, if for every $[\begin{smallmatrix} d_w \\ d_v
\end{smallmatrix}]\in\LTE \times\LTE$ there
exist unique $[\begin{smallmatrix} w \\v \end{smallmatrix}] \in\LTE \times\LTE $ such that
\begin{align}\label{eq:feedbackinterconnection}
  \left\{
  \begin{aligned}
    w &= \Delta v + d_w \\
    v &= G w + d_v 
  \end{aligned}\right.,
\end{align}
and the closed-loop map
$[\![G, \Delta]\!] = ([\begin{smallmatrix} d_w \\ d_v \end{smallmatrix}] \in \LTE \times \LTE) \mapsto (\left[\begin{smallmatrix} w \\ v \end{smallmatrix}\right]
\in \LTE \times\LTE))$
is causal, then the feedback interconnection is called {\em well-posed}.
Moreover, if the induced norm of the restriction to $\LT$ is also
bounded (i.e., $\left\|[\![G,\Delta]\!]\right\|<+\infty$), then the closed-loop is said
to be stable.
\begin{proposition}[\cite{MegRan}] \label{prop:IQCstab}
Let $G:\LTE\rightarrow\LTE$ and $\Delta:\LTE\rightarrow\LTE$ be stable operators. Suppose that $[\![G,\alpha \Delta]\!]$ is well-posed and that $\alpha\Delta\in\mathrm{IQC}(\Pi)$  for all $\alpha\in[0,1]$ for the given self-adjoint multiplier $\Pi:\LTE\rightarrow\LTE$. If there exists $\epsilon>0$ such that 
$$\left\langle 
    \begin{bmatrix} G w \\ w
    \end{bmatrix}, 
    \Pi 
    \begin{bmatrix} G w \\  w \end{bmatrix}
    \right\rangle_{\LT} \leq -\epsilon\|w\|_{\LT}^2 ~\text{for all}~ w\in\LT,
$$
then $[\![G,\Delta]\!]$ is stable.
\end{proposition}

\section{Asynchronous Sample-and-Hold} \label{sec:asynchSH}

The event sequence
$(t_k)_{k\in\mathbb{N}_0}\subset\mathbb{R}_{\geq 0}$ is admissible provided
$t_0=0$, $t_{k+1}-t_{k}>0$ for $k\in\mathbb{N}_0$, and $\lim_{k\rightarrow +\infty}t_{k}=+\infty$. Given admissible event sequence $T=(t_k)_{k\in\mathbb{N}_0}$, the following notation applies:
\begin{itemize}
\item[(i)] $\SSS_{T}$ denotes the sampling operator that maps the continuous-time signal
  $v \in \mathbf{C}_{r}\cap\LTE$ to the discretely indexed signal
  $\tilde{v} =(\tilde{v}_k)_{k\in\mathbb{N}_0}$, such that
  $\tilde{v}_k = v(t_k)$; and
\item[(ii)] ${\HHH_T}$ denotes the hold operator that maps the discretely indexed signal
  $\tilde{v} = (\tilde{v}_k)_{k\in\mathbb{N}_0}$ to the
  continuous-time signal $v\in\LTE$ such that
  $v(t)=(\HHH_T\tilde{v})(t) = \tilde{v}_k$ for  $t \in [t_k,t_{k+1}),\ k\in\mathbb{N}_0$.
Note that $v\in\LTE$ because every finite truncation of the sequence
$\tilde{v}$ is square summable.
\end{itemize}
Synchronous sample-and-hold operators correspond to any composition of the form $\HHH_T\SSS_T$. It is well known that this is equivalent to a
saw-tooth time-varying delay operator~\cite{frid04}.
\begin{lemma}\label{lem:tvdelay}
 Define
  $m(t) = \max\{k~|~t_k\in[0,t],~k\in\mathbb{N}_0\}$ and
  $\sigma(t)=t-t_{m(t)}$ for $t\in\mathbb{R}_{\geq 0}$. 
  Further, for $v\in\LTE$, let
  \begin{align*}
  (\RRR_{\sigma}v)(t) = \begin{cases} v(t-\sigma(t)) &
    \text{if } t-\sigma(t) \geq 0,\\ 0 & \text{otherwise}. \end{cases}
  \end{align*}
  Then
  $\RRR_{\sigma}y = \HHH_{T}\SSS_{T}y$ for
  $y \in \mathbf{C}_r\cap\LTE$. 
\end{lemma}
\begin{IEEEproof}
  Observe that
  $(\HHH_T\SSS_Ty)(t)=y(t_{m(t)}) = y(t-(t-t_{m(t)}))= (\RRR_\sigma y)(t)$ for
  $t\in\mathbb{R}_{\geq 0}$.
\end{IEEEproof}
\begin{remark} \label{rem:commute}
  Note that $\sigma:\mathbb{R}_{\geq 0}\rightarrow\mathbb{R}_{\geq 0}$
  in Lemma~\ref{lem:tvdelay} is piecewise linear with $\sigma(t_k)=0$ for $t_k\in T$. The
  discontinuities are limited to $T$ and the derivative is $1$ almost everywhere. Further,  
  the time-varying delay
  operator $\RRR_{\sigma}$ commutes with multiplication by any
  constant $K$; i.e., $K\RRR_{\sigma} - \RRR_{\sigma} K=0$.
\end{remark}

\begin{figure}[htb]
\centering
\input{asynch_SAH.tex}
\caption{Asynchronous sample-and-hold by composition: $y^\prime = \HHH_{T^\prime}\SSS_{T^\prime}y$ and $y^\star = \HHH_{T^\star}\SSS_{T^\star}y^\prime$. See Lemma~\ref{lem:compo} for $\sigma^\prime$, $\sigma^{\prime\prime}$, and the proof of Lemma~\ref{iqc:gain} for $\psi_l$.}
\label{fig:asynchSAH}
\end{figure}
Given admissible input sample and output update event sequences $T^\prime=(t_k^\prime)_{k\in\mathbb{N}_0}$ and $T^\star=(t_k^\star)_{k\in\mathbb{N}_0}$, the composition 
$\HHH_{T^\star}\SSS_{T^\star} \HHH_{T^\prime}\SSS_{T^\prime}$ yields asynchronous sample-and-hold behaviour of the kind discussed in the introduction; see Fig.~\ref{fig:asynchSAH}.
In view of Lemma~\ref{lem:tvdelay},
$(\HHH_{T^\star}\SSS_{T^\star})(\HHH_{T^\prime}\SSS_{T^\prime}) y = \RRR_{\sigma^\star} \RRR_{\sigma^\prime}y$ for 
  $y\in\mathbf{C}_{r}\cap\LTE$,
where $\sigma^\star(t) = t-t^\star_{n(t)}$ and
$\sigma^\prime(t) = t - t^\prime_{p(t)}$ for $t\in\mathbb{R}_{\geq 0}$, with
\begin{align}
& n(t) = \max\{ k ~|~t^\star_k\in[0,t],~k\in\mathbb{N}_0\} \label{eq:nnn}
\intertext{and}
& p(t) = \max\{ k ~|~t^\prime_k\in
    [0,t],~k\in\mathbb{N}_0\}.
\label{eq:ppp}
\end{align}
\begin{lemma} \label{lem:compo}
With reference to (\ref{eq:nnn}) and (\ref{eq:ppp}), define 
$q(t) = p(t^\star_{n(t)})$ and $\sigma^{\prime\prime}(t) = t - t^\prime_{q(t)}$ for $t\in\mathbb{R}_{\geq 0}$. Then $\RRR_{\sigma^\star}\RRR_{\sigma'} = \RRR_{\sigma^{\prime\prime}}$.  
\end{lemma}
\begin{IEEEproof}
For $y\in\LTE$, note that $y^\star(t) = (\RRR_{\sigma^\star}\RRR_{\sigma^\prime} y)(t) = y^\prime(t-\sigma^\star(t))$, where $y^\prime(t) = (\RRR_{\sigma^\prime}y)(t)=y(t-\sigma^\prime(t))$ for $t\in\mathbb{R}_{\geq 0}$. That is, $y^\star(t) = y(t-\sigma^\star(t)-\sigma^\prime(t-\sigma^\star(t)))$. Now by definition, \begin{align*}\sigma^\star(t) + \sigma^\prime(t-\sigma^\star(t)) &= t - t^\star_{n(t)} + (t - (t - t^\star_{n(t)}) - t^\prime_{p(t-(t-t^\star_{n(t)}))})\\
&= t - t^\prime_{q(t)} = \sigma^{\prime\prime}(t).
\end{align*}
As such, $y^\star(t) = y(t-\sigma^{\prime\prime}(t)) = (\RRR_{\sigma^{\prime\prime}}y)(t)$.
\end{IEEEproof}

\begin{remark}
  Note that
  $\sigma^{\prime\prime}:\mathbb{R}_{\geq
    0}\rightarrow\mathbb{R}_{\geq 0}$ in Lemma~\ref{lem:compo} is also
  piecewise linear, with discontinuities limited to $T^\star$. Moreover, the derivative is $1$ almost everywhere. But $\sigma^{\prime\prime}(t_k^\star)$ may be non-zero and varying for $t_k^\star\in T^\star$; see Fig.~\ref{fig:asynchSAH}. By constrast $\sigma^\prime(t^\prime_k)=0$ for $k\in\mathbb{N}_0$.
\end{remark}

By way of example, consider bounded time-varying delay in a communication link between a sample input buffer and a remote zero-order-hold output buffer. This can be modelled by the asynchronous sample-and-hold composition
$\HHH_{T^\star}\SSS_{T^\star}\HHH_{T^\prime}\SSS_{T^\prime}$ with
$T^\star = \{ t_k^\star = t_k^\prime+\delta_k \}_{k\in\mathbb{N}_0}$
for the given input sample event sequence
$T^\prime=(t^\prime_k)_{k\in\mathbb{N}_0}$, provided the delay satisfies
$0 \leq \delta_k < t^\prime_{k+1} - t^\prime_k$ for every $k\in\mathbb{N}_0$, so that $T^\star$ is an admissible event sequence (e.g., the minimum sample interval is greater than the maximum delay). 
As another example, the asynchronous sample-and-hold operator $\HHH_{T^\star}\SSS_{T^\star}\HHH_{T^\prime}\SSS_{T^\prime}$ could be used to model a shared buffer that is (over-)written at each (possibly aperiodic) time in $T^\prime$, and read asynchronously at times in $T^\star$ to update the output (e.g., in response to write events by the execution of a mechanism that takes a variable time to complete, possibly exceeding the subsequent sample interval).

Of course, the composition of two different synchronous sampled-and-hold operators is just one possible way to arrive at the piecewise linear delay features of asynchronous sample-and-hold operators (e.g., see~\cite{FabbroThesis} where a structured model is proposed that accommodates sample re-ordering by virtue of delay that can exceed the minimum sample interval). Nonetheless, given the preceding examples of practical relevance, attention is restricted here to the class arising via the composition $\HHH_{T^\star}\SSS_{T^\star}\HHH_{T^\prime}\SSS_{T^\prime}$. Exact description of such an operator requires knowledge of the input sample event sequence $T^\prime$ and the output hold update event sequence $T^\star$. These are typically not available ahead of time. On the other hand, it is often possible to bound the inter-sample interval, inter-update interval, and other relationships between the corresponding events. 

The subsequent developments pertain to uncertain asynchronous
sample-and-hold operators defined by the composition
$\HHH_{T^\star}\SSS_{T^\star}\HHH_{T^\prime}\SSS_{T^\prime}$ and the
tuple
$(\tau^\prime,\tau^\star,\tau^\circ,\tau^\natural)\in\mathbb{R}_{>0}
\times \mathbb{R}_{>0}
\times[0,\tau^\star]\times[0,\min\{\tau^\circ,\tau^\prime\}]$ of
bounds such that the event sequences
$T^\prime=(t^\prime_k)_{k\in\mathbb{N}_0}$ and  $T^\star=
(t^\star_k)_{k\in\mathbb{N}_0}$ satisfy the following constraints,
where $\phi_k= \min \{t ~|~ t\in T^\star \cap [t'_k, +\infty)\}$ and $k\in\mathbb{N}_0$:
\begin{align}\label{eq:bounds}
  \begin{split}
    & 0 < t_{k+1}^\prime - t_k^\prime \leq \tau^\prime, \quad
    0 < t_{k+1}^\star - t_k^\star \leq \tau^\star,  \\
    & 
    \phi_k-t_k^\prime \leq \tau^\circ
    \leq \tau^\star, \text{ and } \\
    & \phi_k - \max\{ t ~|~ t\in T'\cap [0,\phi_k]\} \leq \tau^\natural \le \min\{\tau^\circ,\tau^\prime\}.
  \end{split}
\end{align}
The bounds $\tau^\prime$ and $\tau^\star$ are limits on inter-sample and inter-update intervals, respectively. The bound $\tau^\circ$ limits the interval between any sample event and the subsequent hold update event, whereas $\tau^\natural$ limits the interval between such updates and the preceding sample event. As such, these two bounds reflect the degree of asychrony between the sample and hold event sequences. The bound $\tau^\natural$ is not considered in~\cite{CanFabKao18,CanFabKao19}. In general, $\tau^\natural\leq \tau^\circ$, since there could be multiple sample events per hold update (and vice versa). 
\begin{remark}
With reference to Fig.~\ref{fig:asynchSAH}, it follows that $\tau^\prime+\tau^\circ$ is an upper bound on the intervals of time between resets of the time-varying delay $\sigma^{\prime\prime}$, and $\tau^\natural$ is an upper bound on the value to which it resets. 
See the proof of Lemmas~\ref{iqc:gain} and~\ref{iqc:pass} in the next section, where both bounds play an important role in the derivation of IQCs for the class of asynchronous sample-and-hold operators considered.
\end{remark}

Let $F:\LTE\rightarrow\LTE$ be any stable low-pass linear time-invariant
system (i.e., $F$ has strictly proper
transfer function with no poles in the closed right-half plane). Then on $\LTE$,
$\HHH_{T^\star}\SSS_{T^\star}\HHH_{T^\prime}\SSS_{T^\prime} F =
(\mathscr{I}-\Delta) \mathscr{I}^\dagger F$,
where 
\begin{align}
\label{eq:delta}
  \Delta =
(\mathsf{id}-\HHH_{T^\star}\SSS_{T^\star}\HHH_{T^\prime}\SSS_{T^\prime})
  \mathscr{I} = (\mathsf{id} - \RRR_{\sigma^{\prime\prime}})\mathscr{I},
\end{align}
$\RRR_{\sigma^{\prime\prime}}$ is defined in Lemma~\ref{lem:compo} for the realizations $T^\prime$ and $T^\star$ of the sample and hold event sequences, and
$\mathsf{id}$ is the identity on $\LTE$. The operators
$\mathscr{I}:\LTE\rightarrow\LTE$ and
$\mathscr{I}^{\dagger}:\mathcal{D}\rightarrow\LTE$ denote integration
and differentiation, respectivley; i.e.,
$\mathscr{I}=(v\in\LTE) \mapsto ( (t \mapsto \int_0^t v(x)\dif
x)\in\LTE)$ and $\mathscr{I}\mathscr{I}^{\dagger}y = y$ for all $y$ in
the subspace $\mathcal{D}\subset\LTE$ of piecewise differentiable functions.
Note that the range of $F$ is contained in $\mathcal{D}$ and that $\mathscr{I}^\dagger F$ is stable with proper transfer function $sF(s)$. The properties of $F$ also ensure boundedness of the input sampling
operation $\SSS_{T^\prime}F$ on the space $\LT$~\cite{ChenFrancis}. Finally, note that $\Delta$ is causal.

\section{IQC Based Covers (Main Results)} \label{sec:IQCcovers}

 Bounded gain and input feedforward passivity properties are now established for the uncertain operator $\Delta$ in (\ref{eq:delta}), given the tuple $(\tau^\prime,\tau^\star,\tau^\circ,\tau^\natural)\in\mathbb{R}_{>0} \times \mathbb{R}_{>0} \times[0,\tau^\star]\times[0,\min\{\tau^\circ,\tau^\prime\}]$ of bounds such that (\ref{eq:bounds}) holds for the corresponding sample and hold event sequences $T^\prime$ and $T^\star$. Combining these leads to the main result, in which a family of IQCs is provided for $\Delta$.

\begin{lemma}[Bounded Gain]
\label{iqc:gain}
For every $v\in\LT$,
$$\|\Delta\, v\|_{\LT} \leq \left(2(\tau^\prime + \tau^\circ)/\pi + \sqrt{
  (\tau^\prime + \tau^\circ)
  \tau^\natural} \right)\|v\|_{\LT}.$$
\end{lemma}
\begin{IEEEproof}
  Let $\theta_0=0=t_0^\prime=t_0^\star$. With $l=0$ even, define
  \begin{align*}
    \theta_{l+1} &=
    \begin{cases}
      \min \{t ~|~ t\in T^\prime \cap (\theta_l,\infty)\} & \text{if } l \text{ is even}\\
      \min \{ t~|~t\in T^\star \cap [\theta_l,\infty)\} & \text{if } l \text{ is odd}
    \end{cases},~ \quad \psi_l = \theta_{2l},
  \end{align*}
  and $\lambda_l=\max T^\prime \cap [0,\psi_l]$ for
  $l\in\mathbb{N}_0$. Each $\psi_l$ corresponds to a hold output update with the most recent sample, taken at $\lambda_l$; see Fig.~\ref{fig:asynchSAH}. Further, $\psi_{l-1} < \lambda_l \leq \psi_l$ for $l\in\mathbb{N}$, and $\psi_l\rightarrow+\infty$
  as $l\rightarrow+\infty$. 
  
  Let
  $w = \Delta v$. 
  With reference to (\ref{eq:delta}) and Lemma~\ref{lem:compo}, note that  
  $w(t) 
  = \int_{\lambda_l}^tv(x)\dif x$ for $t\in [\psi_l,\psi_{l+1})$ and $l\in\mathbb{N}_0$. So with 
  \begin{align*}
  \underline{w}_l(t) = \int_{\psi_l}^{t} v(x) \dif x
  \quad \text{and} \quad
  \overline{w}_l(t) = \int_{\lambda_l}^{\psi_l} v(x) \dif x
  \end{align*}
  for $t\in [\psi_l,\psi_{l+1})$, and $\underline{w}_l(t)=\overline{w}_l(t)=0$ otherwise, it follows that $w=\overline{w}+\underline{w}$,
  where $\overline{w} = \sum_{l\in\mathbb{N}_0} \overline{w}_l$ and
  $\underline{w}=\sum_{l\in\mathbb{N}_0} \underline{w}_l$.
  Since $\underline{w}(\psi_l)=0$, and
  $\frac{\dif\,\,}{\dif\, t}{\underline{w}}(t)=v(t)$ for $t\in(\psi_l,\psi_{l+1})$, it follows as also noted in~\cite[Lem.~3.2]{Liu10} that application of Wirtinger's inequality~\cite[Thm.~256]{Har52}, with $v\in\LT$, gives
  \begin{align}
    \|\underline{w}\|_{\LT}^2 &= \sum_{l\in\mathbb{N}_0}  \int_{\psi_l}^{\psi_{l+1}} \underline{w}(x)^2 \dif x \nonumber \\ 
    &\leq \sum_{l\in\mathbb{N}_0} \frac{4(\psi_{l+1}-\psi_l)^2}{\pi^2} \int_{\psi_{l}}^{\psi_{l+1}} 
    v(x)^2\dif x \nonumber \\
    &\leq\left(\frac{2(\tau^\prime+\tau^\circ)}{\pi}
    \|v\|_{\LT}\right)^2. \label{eq:wul_bound}
  \end{align}
  The last inequality above holds because
  $(\psi_{l+1}-\psi_l) 
  =
  (\theta_{2l+1}-\theta_{2l}) + (\theta_{2(l+1)}-\theta_{2l+1}) \leq
  \tau^\prime + \tau^\circ$ for
  $l\in\mathbb{N}_0$ by (\ref{eq:bounds}). Further,
  \begin{align}
    \|\overline{w}\|_{\LT}^2 &= \sum_{l\in\mathbb{N}_0} (\psi_{l+1} -\psi_l) (\int_{\lambda_l}^{\psi_l} v(x) \dif x)^2 \nonumber \\ 
    &\leq \sum_{l\in\mathbb{N}_0} (\psi_{l+1}-\psi_l) (\psi_l-\lambda_l)\int_{\lambda_l}^{\psi_l} v(x)^2 \dif x \nonumber \\
    & \leq \left(\sqrt{(\tau^\prime+\tau^\circ) \tau^\natural} \cdot  \|v\|_{\LT}\right)^2. \label{eq:wol_bound}
  \end{align}
  The equality above holds because $\overline{w}(t)$ is constant for $t\in[\psi_l,\psi_{l+1})$. The first inequality holds by application of 
  Jensen's inequality~\cite[Thm.~3.3]{Rudin}, and the second because
  $(\psi_l-\lambda_l)\leq \tau^\natural$ by (\ref{eq:bounds}) and
  $\int_{\lambda_l}^{\psi_l} v(x)^2 \dif x \leq
  \int_{\lambda_l}^{\lambda_{l+1}} v(x)^2 \dif x$.
  
  In view of (\ref{eq:wul_bound}) and (\ref{eq:wol_bound}), it follows that $\underline{w},\overline{w}\in\LT$, and hence, $w=\underline{w}+\overline{w}\in\LT$. Further, the claimed gain bound holds by the triangle inequality; i.e., $\|w\|_{\LT} \leq \|\underline{w}\|_{\LT}
  + \|\overline{w}\|_{\LT}$.
\end{IEEEproof}

\begin{remark}
  Synchronous sample-and-hold corresponds to $T^\prime = T^\star$, and thus, $\tau^\circ=\tau^\natural=0$. In this case, the $\LT$ gain bound in
  Theorem~\ref{iqc:gain} is $2\tau^\prime/\pi$, which is exactly the induced
  norm of the corresponding $\Delta$, as shown in~\cite{Mir07}; i.e., the bound is tight. The gain bound provided in the preliminary work~\cite{CanFabKao18,CanFabKao19} corresponds to taking $\tau^\natural=\tau^\circ$. However, this can be conservative. For example, if the hold output is simply a down-sampled version of the input, then $\tau^\natural$ can be taken to be zero, although $\tau^\circ$ is non-zero, making the gain bound tighter.
\end{remark}

\begin{lemma}[Input Feedforward Passivity]
\label{iqc:pass}
For $v\in\LT$, the following holds:
$~\langle \Delta v,v \rangle_{\LT} + (\tau^\natural /2)\|v\|_{\LT}^2 \ge 0.$
\end{lemma}
\begin{IEEEproof}
  Using notation from the proof of
  Lemma~\ref{iqc:gain}, with $w=\Delta v\in\LT$, $\frac{\dif\,\,}{\dif t}{w}(t) = v(t)$, and thus,
  $w(t)v(t)=\frac{1}{2}\frac{\dif\,\,}{\dif t}{z}(t)$, where $z=w^2$, for $l\in\mathbb{N}_0$ and
  $t\in(\psi_l,\psi_{l+1})$. As such, 
  $$\int_{\psi_l}^{\psi_{l+1}} w(x)v(x)\dif x \\
  = \frac{1}{2} \left(
    (\int_{\lambda_l}^{\psi_{l+1}} v(x)
    \dif x )^2 - (\int_{\lambda_l}^{\psi_l} v(x)\dif
    x )^2\right).$$
    This implies
    $$0\leq \int_{\psi_l}^{\psi_{l+1}} w(x)v(x)\dif x +
    \frac{1}{2}(\int_{\lambda_l}^{\psi_l} v(x)\dif x )^2.$$ 
    By application of Jensen's inequality~\cite[Thm.~3.3]{Rudin}, it follows
    that
    \begin{align*}
    0 &\leq \int_{\psi_l}^{\psi_{l+1}} w(x)v(x)\dif x +
    ((\psi_l-\lambda_l)/2)\int_{\lambda_l}^{\psi_l} v(x)^2\dif x \\
    &\leq
    \int_{\psi_l}^{\psi_{l+1}} w(x)v(x)\dif x +
    ((\psi_l-\lambda_l)/2)\int_{\lambda_l}^{\lambda_{l+1}} v(x)^2\dif
    x.
    \end{align*}
As $(\psi_l-\lambda_l)\leq \tau^\natural$ by
(\ref{eq:bounds}), summing over $l\!\in\!\mathbb{N}_0$ yields the result.
\end{IEEEproof}
\begin{remark}
Lemma~\ref{iqc:pass} reveals that $(\tau^\natural/2) \mathsf{id} + \Delta$ is passive; i.e., $-(\tau^\natural/2)\leq 0$ is a lower bound for the input feedforward passivity index~\cite{Rahn18}. This is tighter than the one that follows from~\cite[Lem.~5]{CanFabKao18}, where $\tau^\natural=\tau^\circ$, giving the weaker constraint $\langle \Delta v,v \rangle_{\LT} \geq -(\tau^\circ/2)\|v\|_{\LT}^2$, since $\tau^\natural \leq \tau^\circ$ in general.
\end{remark}

\begin{theorem}
\label{thm:IQCsp}
Let
$\beta = (~ 2(\tau^\prime+\tau^\circ) / \pi + \sqrt{
  (\tau^\prime+\tau^\circ) \tau^\natural} ~)^2$ and 
\begin{align*}
  \Pi(X,Y)= &\left( \begin{bmatrix} v \\ w \end{bmatrix} \! \in \! \LT
     \mapsto  ( t\! \in \! \mathbb{R}_{\geq 0} \mapsto 
\begin{bmatrix}
    \beta X + \tau^\natural Y & Y \\ Y & -X
  \end{bmatrix}\! \begin{bmatrix} v(t) \\ w(t) \end{bmatrix}) \! \in\! \LT \!\! \right)\!\!.
\end{align*}
Then $\Delta\in \mathrm{IQC}(\Pi(X,Y))$
for every $X\geq 0$ and $Y\ge 0$. 
\end{theorem}
\begin{IEEEproof}
  For arbitrary $v\in\LT$, let $w=\Delta v$. Note that
\begin{align}
\begin{split}\label{eq:bigiqc}
&\left \langle \begin{bmatrix} v \\ w \end{bmatrix}, \Pi(X,Y) 
  \begin{bmatrix} v \\ w \end{bmatrix}
\right\rangle_{\LT}
= \beta \|X^{1/2}v\|_{\LT}^2-\|X^{1/2}w\|_{\LT}^2
+2 \langle Y^{1/2}w, Y^{1/2}v\rangle_{\LT}
+ \tau^\natural \norm{Y^{1/2}v}_{\LT}^2.
\end{split}
\end{align}
Since the time-varying delay $\RRR_{\sigma^{\prime\prime}}$,
and thus 
$\Delta=(\mathrm{Id}-\RRR_{\sigma^{\prime\prime}})\mathscr{I}$, 
commute with multiplication by a constant gain (see Remark~\ref{rem:commute}), $X^{1/2} w = \Delta X^{1/2} v$ and
$Y^{1/2} w = \Delta Y^{1/2} v$. With Lemma~\ref{iqc:gain}, this
implies
$\beta \|X^{1/2}v\|_{\LT}^2-\|X^{1/2}w\|_{\LT}^2\geq 0$  in (\ref{eq:bigiqc}).
Similarly, the remaining terms in (\ref{eq:bigiqc}) are non-negative by Lemma~\ref{iqc:pass}. In summary,
$\left\langle [\begin{smallmatrix} v \\ \Delta v \end{smallmatrix}],
  \Pi(X,Y) [\begin{smallmatrix} v \\ \Delta v \end{smallmatrix}]
\right\rangle_{\LT} \geq 0$
for $v\in\LT$. 
\end{IEEEproof}

\begin{figure}[htbp]
\centering
\setlength{\unitlength}{1500sp}%
\begingroup\makeatletter\ifx\SetFigFont\undefined%
\gdef\SetFigFont#1#2#3#4#5{%
  \reset@font\fontsize{#1}{#2pt}%
  \fontfamily{#3}\fontseries{#4}\fontshape{#5}%
  \selectfont}%
\fi\endgroup%
\begin{picture}(10074,8772)(139,-8773)
\put(1501,-4536){\makebox(0,0)[b]{\smash{{\SetFigFont{11}{16.8}{\rmdefault}{\mddefault}{\updefault}{\color[rgb]{0,0,0}$W$}%
}}}}
\put(1501,-6586){\makebox(0,0)[b]{\smash{{\SetFigFont{11}{16.8}{\rmdefault}{\mddefault}{\updefault}{\color[rgb]{0,0,0}$\mathscr{I}^\dagger F$}%
}}}}
{\color[rgb]{0,0,0}\thinlines
\put(4951,-961){\circle{300}}
}%
{\color[rgb]{0,0,0}\put(2701,-1561){\framebox(1200,1200){}}
}%
{\color[rgb]{0,0,0}\put(4801,-961){\vector(-1, 0){900}}
}%
{\color[rgb]{0,0,0}\put(6001,-961){\vector(-1, 0){900}}
}%
{\color[rgb]{0,0,0}\put(4951,-61){\vector( 0,-1){750}}
}%
{\color[rgb]{0,0,0}\put(2251,-961){\line( 0,-1){1650}}
\put(2251,-2611){\vector( 1, 0){1200}}
}%
{\color[rgb]{0,0,0}\put(7651,-2611){\line( 1, 0){1050}}
\put(8701,-2611){\line( 0, 1){1650}}
\put(8701,-961){\line(-1, 0){2700}}
}%
{\color[rgb]{0,0,0}\put(3451,-3211){\framebox(4200,1200){}}
}%
{\color[rgb]{0,0,0}\put(2701,-961){\line(-1, 0){450}}
}%
\put(5551,-2686){\makebox(0,0)[b]{\smash{{\SetFigFont{11}{16.8}{\rmdefault}{\mddefault}{\updefault}{\color[rgb]{0,0,0}$\HHH_{T^\star}\SSS_{T^\star}
          \HHH_{T^\prime}\SSS_{T^\prime}F$}%
}}}}
\put(5176,-211){\makebox(0,0)[b]{\smash{{\SetFigFont{11}{16.8}{\rmdefault}{\mddefault}{\updefault}{\color[rgb]{0,0,0}$d$}%
}}}}
\put(3301,-1066){\makebox(0,0)[b]{\smash{{\SetFigFont{11}{16.8}{\rmdefault}{\mddefault}{\updefault}{\color[rgb]{0,0,0}$P$}%
}}}}
\put(4726,-736){\makebox(0,0)[b]{\smash{{\SetFigFont{11}{16.8}{\rmdefault}{\mddefault}{\updefault}{\color[rgb]{0,0,0}+}%
}}}}
\put(4351,-861){\makebox(0,0)[b]{\smash{{\SetFigFont{11}{16.8}{\rmdefault}{\mddefault}{\updefault}{\color[rgb]{0,0,0}$u$}%
}}}}
\put(2026,-1111){\makebox(0,0)[b]{\smash{{\SetFigFont{11}{16.8}{\rmdefault}{\mddefault}{\updefault}{\color[rgb]{0,0,0}$y$}%
}}}}
\put(5251,-1111){\makebox(0,0)[b]{\smash{{\SetFigFont{11}{16.8}{\rmdefault}{\mddefault}{\updefault}{\color[rgb]{0,0,0}\_}%
}}}}
{\color[rgb]{0,0,0}\put(5401,-5161){\circle{300}}
}%
{\color[rgb]{0,0,0}\put(9151,-6511){\circle{300}}
}%
{\color[rgb]{0,0,0}\put(5251,-5161){\vector(-1, 0){900}}
}%
{\color[rgb]{0,0,0}\put(6451,-5161){\vector(-1, 0){900}}
}%
{\color[rgb]{0,0,0}\put(3151,-5161){\line(-1, 0){450}}
}%
{\color[rgb]{0,0,0}\put(2101,-6511){\vector(-1, 0){  0}}
\put(2101,-6511){\vector( 1, 0){2475}}
}%
{\color[rgb]{0,0,0}\put(9151,-6361){\line( 0, 1){1200}}
\put(9151,-5161){\line(-1, 0){2700}}
}%
{\color[rgb]{0,0,0}\put(5401,-4411){\vector( 0,-1){600}}
}%
{\color[rgb]{0,0,0}\put(5401,-4411){\line( 1, 0){4800}}
}%
{\color[rgb]{0,0,0}\put(3151,-5761){\framebox(1200,1200){}}
}%
{\color[rgb]{0,0,0}\put(301,-6511){\line( 0,-1){1650}}
\put(301,-8161){\vector( 1, 0){3900}}
}%
{\color[rgb]{0,0,0}\put(6151,-8161){\line( 1, 0){3900}}
\put(10051,-8161){\line( 0, 1){1650}}
\put(10051,-6511){\vector(-1, 0){750}}
}%
{\color[rgb]{0,0,0}\put(301,-6511){\line( 1, 0){600}}
}%
{\color[rgb]{0,0,0}\put(4201,-8761){\framebox(1950,1200){}}
}%
{\color[rgb]{0,0,0}\put(5776,-6511){\vector( 1, 0){3225}}
}%
{\color[rgb]{0,0,0}\put(4576,-7111){\framebox(1200,1200){}}
}%
{\color[rgb]{0,0,0}\put(2701,-5161){\line( 0,-1){1350}}
}%
{\color[rgb]{0,0,0}\put(4876,-5161){\line( 0, 1){750}}
\put(4876,-4411){\vector(-1, 0){2775}}
}%
{\color[rgb]{0,0,0}\put(901,-5011){\framebox(1200,1200){}}
}%
\thicklines
{\color[rgb]{0.7,0.7,0.7}\put(706,-7156){\oval(210,210)[bl]}
\put(706,-3766){\oval(210,210)[tl]}
\put(9646,-7156){\oval(210,210)[br]}
\put(9646,-3766){\oval(210,210)[tr]}
\put(706,-7261){\line( 1, 0){8940}}
\put(706,-3661){\line( 1, 0){8940}}
\put(601,-7156){\line( 0, 1){3390}}
\put(9751,-7156){\line( 0, 1){3390}}
}%
\thinlines
{\color[rgb]{0,0,0}\put(901,-4411){\vector(-1, 0){750}}
}%
\put(10051,-4261){\makebox(0,0)[b]{\smash{{\SetFigFont{11}{16.8}{\rmdefault}{\mddefault}{\updefault}{\color[rgb]{0,0,0}$d$}%
}}}}
\put(5176,-8236){\makebox(0,0)[b]{\smash{{\SetFigFont{11}{16.8}{\rmdefault}{\mddefault}{\updefault}{\color[rgb]{0,0,0}$\Delta$}%
}}}}
\put(301,-6361){\makebox(0,0)[b]{\smash{{\SetFigFont{11}{16.8}{\rmdefault}{\mddefault}{\updefault}{\color[rgb]{0,0,0}$v$}%
}}}}
\put(10051,-6361){\makebox(0,0)[b]{\smash{{\SetFigFont{11}{16.8}{\rmdefault}{\mddefault}{\updefault}{\color[rgb]{0,0,0}$w$}%
}}}}
\put(3751,-5266){\makebox(0,0)[b]{\smash{{\SetFigFont{11}{16.8}{\rmdefault}{\mddefault}{\updefault}{\color[rgb]{0,0,0}$P$}%
}}}}
\put(5176,-4936){\makebox(0,0)[b]{\smash{{\SetFigFont{11}{16.8}{\rmdefault}{\mddefault}{\updefault}{\color[rgb]{0,0,0}+}%
}}}}
\put(5176,-6586){\makebox(0,0)[b]{\smash{{\SetFigFont{11}{16.8}{\rmdefault}{\mddefault}{\updefault}{\color[rgb]{0,0,0}$F$}%
}}}}
\put(301,-4261){\makebox(0,0)[b]{\smash{{\SetFigFont{11}{16.8}{\rmdefault}{\mddefault}{\updefault}{\color[rgb]{0,0,0}$z$}%
}}}}
\put(2476,-5311){\makebox(0,0)[b]{\smash{{\SetFigFont{11}{16.8}{\rmdefault}{\mddefault}{\updefault}{\color[rgb]{0,0,0}$y$}%
}}}}
\put(5712,-5284){\makebox(0,0)[b]{\smash{{\SetFigFont{11}{16.8}{\rmdefault}{\mddefault}{\updefault}{\color[rgb]{0,0,0}\_}%
}}}}
\put(9451,-6334){\makebox(0,0)[b]{\smash{{\SetFigFont{11}{16.8}{\rmdefault}{\mddefault}{\updefault}{\color[rgb]{0,0,0}\_}%
}}}}
\put(8862,-6433){\makebox(0,0)[b]{\smash{{\SetFigFont{11}{16.8}{\rmdefault}{\mddefault}{\updefault}{\color[rgb]{0,0,0}+}%
}}}}
{\color[rgb]{0,0,0}\put(901,-7111){\framebox(1200,1200){}}
}%
\end{picture}%
\caption{A sampled-data implementation of a feedback loop (top) and a loop-transformation for robust stability and performance analysis via Proposition~\ref{prop:IQCstab} and Theorem~\ref{thm:IQCsp} (bottom).}
\label{fig:feedbackloop}
\end{figure}
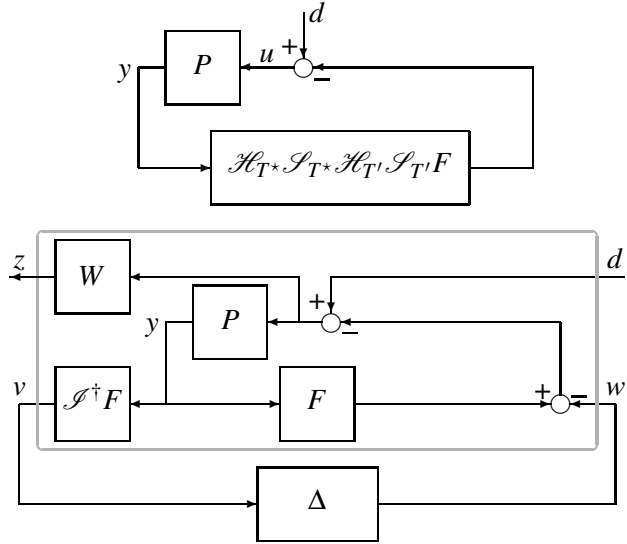

\section{Robust Performance Analysis of an Example Sampled-Data Feedback Loop} \label{sec:examples}

Consider the top part of Fig.~\ref{fig:feedbackloop}, which shows the digital implementation of a disturbance matching feedback loop. The path from the sensor to the actuator comprises a low-pass filter $F$ and a sample buffer of depth one; i.e., it can only store one sample. This buffer is (over-)written with sensor samples by an external trigger. It is read to update the actuator in response to such events with variable delay. This feedback path can be modelled by the asynchronous sampled-and-hold operator $\HHH_{T^\star}\SSS_{T^\star}\HHH_{T^\prime}\SSS_{T^\prime}F$, where the sensor sampling events $T^\prime$ and actuator update events $T^\star$ satisfy \eqref{eq:bounds} for an appropriate tuple of bounds $(\tau^\prime, \tau^\star,\tau^\circ,\tau^\natural)\in\mathbb{R}_{> 0}\times\mathbb{R}_{> 0}\times[0,\tau^\star]\times[0,\min\{\tau^\circ,\tau^\prime\}]$. To illustrate an application of the results in Section~\ref{sec:IQCcovers}, closed-loop stability and performance are considered below with
\begin{align} \label{eq:tauexamp} 
\tau^\prime=h,~ \tau^\star=(1+\delta) h,~ \tau^\circ=\delta h,~\text{and}~ \tau^\natural=h \cdot \min\{\delta,1\},
\end{align}
for $h>0$ and $\delta\geq 0$.
This corresponds to a bound $\delta h\geq 0$ on the actuator update response time, which may exceed the maximum sample interval when $\delta > 1$. In the following, the anti-aliasing filter $F$ and transfer function $P$ are such that the feedback interconnection $[\![P,-F]\!]$ is stable; i.e., the ideal closed-loop with infinitely frequent sampling and synchronous hold is stable. For the sake of argument, the performance weight transfer function in the loop-transformed system at the bottom of Fig.~\ref{fig:feedbackloop} is taken to be $W(s)=F(s)$.


Using \eqref{eq:delta} and Lemma~\ref{iqc:gain}, note that $\HHH_{T^\star}\SSS_{T^\star}
\HHH_{T^\prime} \SSS_{T^\prime}F = F - \Delta \mathscr{I}^\dagger F$ is stable. As such, stability of the feedback interconnection $[\![P,-\HHH_{T^\star}\SSS_{T^\star}
\HHH_{T^\prime} \SSS_{T^\prime}F ]\!]$ in the top of Fig.~\ref{fig:feedbackloop} is equivalent to stability of the closed-loop map $d\mapsto y$. Specifically, any $\LT$ disturbance at the plant output can as such be transferred to the plant input, and if $y$ is an element of $\LT$ for $d\in\LT$, then so is $u=d- \HHH_{T^\star}\SSS_{T^\star}
\HHH_{T^\prime} \SSS_{T^\prime}F y$. Also observe that the feedback interconnection $[\![P,-F]\!]$ is stable, and thus, that if $w\in\LT$ in the loop-transformed system at the bottom of Fig.~\ref{fig:feedbackloop}, then $y\in\LT$. Therefore, verification that the interconnection $[\![( \mathscr{I}^\dagger F \left[ \begin{smallmatrix} 0 & \mathsf{id} \end{smallmatrix} \right]  [\![P,-F]\!] 
\left[ \begin{smallmatrix}  \mathsf{id} \\ 0 \end{smallmatrix} \right])
, \Delta ]\!]$ is stable, implies $[\![P,-\HHH_{T^\star}\SSS_{T^\star}
\HHH_{T^\prime} \SSS_{T^\prime}F ]\!]$ is stable.

Define 
\begin{align*}
G= \begin{bmatrix}
G_{zd} & G_{zw} \\ G_{vd} & G_{vw}
\end{bmatrix}
=
\begin{bmatrix}
W & 0 \\ 0 & \mathscr{I}^\dagger F 
\end{bmatrix}
[\![ P, -F ]\!]
\begin{bmatrix}
\mathsf{id} & \mathsf{id} \\
0 & 0 
\end{bmatrix},
\end{align*}
noting that $G_{zd}$, $G_{zw}$, $G_{vd}$ and $G_{vw}$ are all stable linear time-invariant systems with proper transfer functions. Then $[\![P,-\HHH_{T^\star}\SSS_{T^\star}
\HHH_{T^\prime} \SSS_{T^\prime}F ]\!]$ is stable if $[\![G_{vw},\Delta]\!]$ is stable. Further, the weighted closed-loop performance map $d \mapsto z=Wu$ is given by $G_{zd} + G_{zw} \left[ \begin{smallmatrix} 0 & \mathsf{id} \end{smallmatrix} \right]  [\![G_{vw},\Delta]\!] 
\left[ \begin{smallmatrix}  \mathsf{id} \\ 0 \end{smallmatrix} \right] G_{vd}$, which is bounded in this case. Given any state-space realization
\begin{align*}
G(s) = \begin{bmatrix} C_z \\ C_v \end{bmatrix}(sI-A)^{-1} 
\begin{bmatrix} 
B_d & B_w
\end{bmatrix} + \begin{bmatrix}
D_{zd} & D_{zw} \\ D_{vd} & D_{vw}
\end{bmatrix},
\end{align*}
with $A$ Hurwitz, the following stability and performance certificates are now established via Theorem~\ref{thm:IQCsp}, Proposition~\ref{prop:IQCstab} and the so-called Kalman-Yakubovic-Popov (KYP) lemma~\cite{Ran96}. The following notation applies below: For symmetric matrix $M=M^\top\in\mathbb{R}^{n\times n}$, $M\succ 0$ means there exist $c>0$ such that $x^\top M x \geq c ~x^\top x$ for all $x\in\mathbb{R}^n$, and $M\prec 0$ means $-M\succ 0$.
\begin{theorem}[Robust Stability] \label{thm:robstab}
Given $h>0$ and $\delta\geq 0$, let 
$\beta = (\,2(1+\delta)h / \pi + \sqrt{
  (1+\delta)h^2 \cdot \min\{\delta,1\}}\,)^2$ and $\eta = h\cdot\min\{\delta,1\}$.
If $X\geq 0$, $Y\geq 0$ and $Q=Q^\top$ exist such that 
\begin{align} \label{eq:stabLMI}
&\begin{bmatrix}
A & B_{w}\\
I & 0
\end{bmatrix}^\top
\begin{bmatrix} 
0 & Q \\ Q & 0
\end{bmatrix} 
\begin{bmatrix}
A & B_{w}\\
I & 0
\end{bmatrix}
+ 
\begin{bmatrix}
C_{v} & D_{vw}\\
0 & I
\end{bmatrix}^\top
\begin{bmatrix}
    \beta X + \eta Y & Y \\ Y & -X
  \end{bmatrix}
\begin{bmatrix}
C_{v} & D_{vw}\\
0 & I
\end{bmatrix}
\prec 0,
\end{align}
then $[\![P,-\HHH_{T^\star}\SSS_{T^\star}
\HHH_{T^\prime} \SSS_{T^\prime}F ]\!]$ is stable for all event sequences $T^\prime$ and $T^\star$ consistent with \eqref{eq:bounds} and \eqref{eq:tauexamp}.
\end{theorem}
\begin{IEEEproof}
The interconnection $[\![G_{vw},\alpha \Delta]\!]$ is well-posed for $\alpha\in[0,1]$ since the so-called instantaneous gain of $\alpha\Delta = (v\in\LT\mapsto (t\in\mathbb{R}_{\geq 0} \mapsto \alpha \int_{q(t)}^t v(t)\, \dif x\in\LT))$ is zero; 
see~\cite{Willems71}. Further, since $\|\alpha\Delta v \|=\alpha\|\Delta v\|$ and $\langle \alpha \Delta v, v\rangle_{\LT} = \alpha\langle \Delta v, v\rangle$ for $v\in\LT$,  Lemmas~\ref{iqc:gain} and \ref{iqc:pass} hold with $\alpha\Delta$ in place of $\Delta$ for $\alpha\in[0,1]$. Therefore, as in the proof of Theorem~\ref{thm:IQCsp}, it follows that $\alpha\Delta\in\mathrm{IQC}([\begin{smallmatrix} 
\beta X + \eta Y & Y \\ Y & -X
\end{smallmatrix}])$ for every $X\geq 0$, $Y\geq 0$, and $\alpha\in[0,1]$. Now applying Proposition~\ref{prop:IQCstab}, if there exists $X\geq 0$, $Y\geq 0$ and $\epsilon>0$ such that 
\begin{align} \label{eq:IQCstab_cert}
\left\langle 
    \left[\begin{matrix} G_{vw} w \\ w
    \end{matrix}\right], 
    \left[\begin{matrix} 
\beta X + \eta Y & Y \\ Y & -X
\end{matrix}\right] 
    \left[\begin{matrix} G_{vw} w \\  w \end{matrix}\right]
    \right\rangle_{\LT} \leq -\epsilon\|w\|_{\LT}^2
    \end{align}
    for all $w\in\LT$, then $[\![G_{vw},\Delta]\!]$ is stable, which as discussed above, implies $[\![P,-\HHH_{T^\star}\SSS_{T^\star}\HHH_{T^\prime} \SSS_{T^\prime}F ]\!]$ is stable. As such,
    the stated result holds by application of the KYP lemma~\cite{Ran96}, whereby the existence of $X\geq 0$, $Y\geq 0$ and $\epsilon>0$ such that \eqref{eq:IQCstab_cert} holds for all $w\in\LT$ is equivalent to the existence of $X\geq 0$, $Y\geq 0$ and $Q=Q^\top$ such that \eqref{eq:stabLMI} holds.
\end{IEEEproof}

\begin{theorem}[Robust Performance] \label{thm:robperf}
With the notation of Theorem~\ref{thm:robstab}, given $\gamma >0$, if $X\geq 0$, $Y\geq 0$ and $Q=Q^\top$ exist such that 
\begin{align}
\label{eq:perfLMI}
&\begin{bmatrix}
A & \begin{bmatrix} B_d & B_{w} \end{bmatrix} \\
I & 0
\end{bmatrix}^\top
\begin{bmatrix} 
0 & Q \\ Q & 0
\end{bmatrix} 
\begin{bmatrix}
A & \begin{bmatrix} B_d & B_{w} \end{bmatrix}\\
I & 0
\end{bmatrix}
\nonumber \\
&\quad
+ 
\begin{bmatrix}
\begin{bmatrix} C_z \\
C_v\end{bmatrix} 
& 
\begin{bmatrix}
D_{zd} & D_{zw} \\ D_{vd} & D_{vw}
\end{bmatrix}
\\
0 & I
\end{bmatrix}^\top
\begin{bmatrix}
I & 0 & 0 & 0 \\
0 &  \beta X + \eta Y & 0 & Y\\
0 & 0 & -\gamma^2 I & 0 \\
0 & Y & 0 &  -X
  \end{bmatrix}
\begin{bmatrix}
\begin{bmatrix} C_z \\
C_v\end{bmatrix} 
& 
\begin{bmatrix}
D_{zd} & D_{zw} \\ D_{vd} & D_{vw}
\end{bmatrix}
\\
0 & I
\end{bmatrix}
\prec 0,
\end{align}
then $[\![P,-\HHH_{T^\star}\SSS_{T^\star}
\HHH_{T^\prime} \SSS_{T^\prime}F ]\!]$ is stable and $\|d \mapsto z=W u\| \leq  \gamma$ for all  $T^\prime$ and $T^\star$ that are consistent with \eqref{eq:bounds} and \eqref{eq:tauexamp}.
\end{theorem}
\begin{IEEEproof}
This performance analysis extension of the stability result in Theorem~\ref{thm:robstab} is standard. A brief proof is provided for completeness. By the KYP lemma~\cite{Ran96}, the existence of $X\geq 0$, $Y\geq 0$ and $Q=Q^\top$ such that (\ref{eq:perfLMI}) holds is equivalent to the existence of $X\geq 0$, $Y\geq 0$ and $\epsilon>0$ such that 
\begin{align*}
\left\langle
\begin{bmatrix}
G \left[\begin{smallmatrix}
d \\ w \end{smallmatrix} \right] \\
\left[\begin{smallmatrix}
d \\ w \end{smallmatrix} \right]
\end{bmatrix},
\left[\begin{smallmatrix}
I & 0 & 0 & 0 \\
0 &  \beta X + \eta Y & 0 & Y\\
0 & 0 & -\gamma^2 I & 0 \\
0 & Y & 0 &  -X
\end{smallmatrix} \right]
\begin{bmatrix}
G \left[\begin{smallmatrix}
d \\ w \end{smallmatrix} \right] \\
\left[\begin{smallmatrix}
d \\ w \end{smallmatrix} \right]
\end{bmatrix}
\right\rangle_{\LT} \leq 
-\epsilon\|\left[\begin{smallmatrix}
d \\ w \end{smallmatrix} \right]\|_2^2
\end{align*}
for all $\left[\begin{smallmatrix}
d \\ w \end{smallmatrix} \right]\in\LT$. With $d=0$, this implies (\ref{eq:IQCstab_cert}) holds for all $w\in\LT$, since $\langle G_{zw} w, G_{zw} w\rangle_{\LT}\geq 0$. Thus, $[\![P,-\HHH_{T^\star}\SSS_{T^\star}
\HHH_{T^\prime} \SSS_{T^\prime}F ]\!]$ is stable, as seen in the proof of Theorem~\ref{thm:robstab}. Moreover, with $\left[\begin{smallmatrix} z \\ v \end{smallmatrix}\right] = G \left[\begin{smallmatrix} d \\ w \end{smallmatrix}\right]$ and $w=\Delta v$ as shown in Fig.~\ref{fig:feedbackloop}, it follows that
\begin{align*}
\langle z, z \rangle_{\LT} - \gamma^2 \langle d, d \rangle_{\LT}
&\leq 
-\epsilon \|\! \left[\begin{smallmatrix} d \\ w \end{smallmatrix}\right]\! \|_{\LT}^2
   -\langle 
\left[\begin{smallmatrix} v \\ \Delta v \end{smallmatrix}\right],
\left[\begin{smallmatrix}
\beta X + \eta Y & Y\\
 Y &  -X
\end{smallmatrix} \right]
\left[\begin{smallmatrix} v \\ \Delta v \end{smallmatrix}\right]
\rangle_{\LT} \leq 0,
\end{align*}
since the term involving $v$ is non-negative by Theorem~\ref{thm:IQCsp}. 
\end{IEEEproof}

\begin{remark}
 Verification of 
 (\ref{eq:stabLMI}) or (\ref{eq:perfLMI}) 
 is a standard question of linear matrix inequality (LMI) feasibility. It can be decided by posing a finite-dimensional semi-definite program. To this end, the CVX package~\cite{cvx} with default solver SDPT3~\cite{sdpt3} is used here. 
\end{remark}

\begin{figure}[htbp]
\centering
\includegraphics[width = 0.45\textwidth]{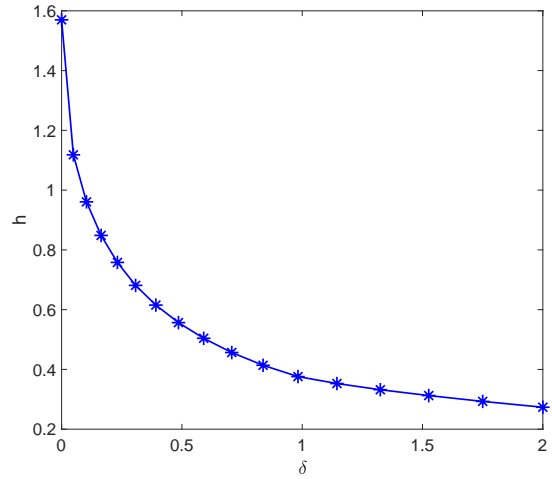} 
\caption{Largest verifiable inter-sample interval bound $h$ for stability as hold update asynchrony bound $\delta$ is varied, for $P(s)=1/s$ and $F(s)=1/(0.1s+1)$.}
\label{fig:example1a}
\end{figure}

\begin{figure}[htbp]
\centering
\includegraphics[width = 0.45\textwidth]{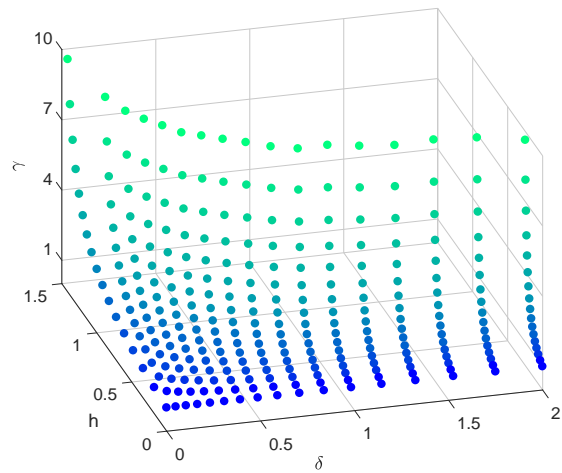} 
\caption{Largest verifiable $\LT$-gain bound $\gamma$ over a grid of stable $(h,\delta)$ pairs for $P(s)=1/s$ and $F(s)=1/(0.1s+1)$.}
\label{fig:example1b}
\end{figure}

Numerical results are given in Figs.~\ref{fig:example1a} and~\ref{fig:example1b} for $P(s)=1/s$ and $F(s)=1/(0.1s+1)$, whereby the nominal interconnection $[\![P,-F]\!]$ is stable.  Fig.~\ref{fig:example1a} shows the largest value of the inter-sample interval bound $h>0$ for which the condition (\ref{eq:stabLMI}) can be verified numerically, as the hold update asynchrony bound $\delta$ is varied from $0$ to $2$. As might be expected, with increasing $\delta$ the largest verifiable $h>0$ decreases. Fig.~\ref{fig:example1b} shows the smallest $\LT$-gain bound $\gamma$ on the closed-loop map $d\mapsto z$, for which the condition (\ref{eq:perfLMI}) can be verified, over a grid of $(h,\delta)$ pairs. As might be expected, the verified gain bound increases sharply as $(h,\delta)$ approaches the approximate ``stability boundary'' shown in Fig.~\ref{fig:example1a}. For this example, it turns out that fixing $Y=0$ in (\ref{eq:stabLMI}) and (\ref{eq:perfLMI}) yields the same results; i.e., only the gain bound part of the IQC from Lemma~\ref{iqc:gain} is important in this example. 

\begin{figure}[htbp]
\centering
\includegraphics[width = 0.45\textwidth]{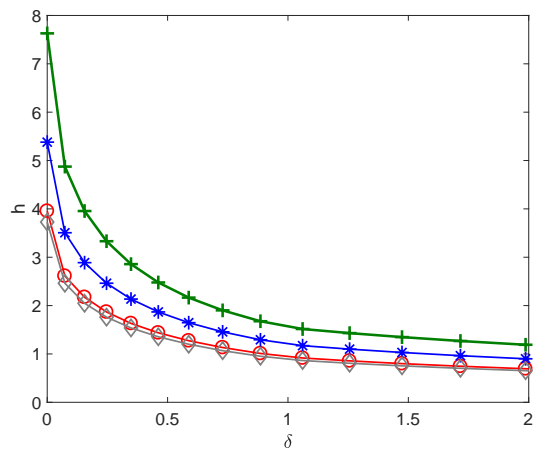} 
\caption{Largest verifiable inter-sample interval bound $h$ for stability as hold update asynchrony bound $\delta$ is varied, for $P(s)=0.9(T_z s-1)/(s^2+2s+1)$ and $F(s)=1/(0.1s+1)$; (blue-star $T_z=0.2$, green-plus $T_z=0.05$) correspond to verification with $X\geq0$ and $Y\geq 0$ in (\ref{eq:stabLMI}); (red-circle $T_z=0.2$, grey-diamond $T_z=0.05$) correspond to verification with $X\geq 0$ and $Y=0$ in (\ref{eq:stabLMI}).}
\label{fig:example2}
\end{figure}

Consider $P(s)=0.9(T_z s-1)/(s^2+2s+1)$ and $F(s)=1/(0.1 s+1)$ with $T_z=0.05$ and $T_z=0.2$, noting that the nominal $[\![P,-F]\!]$ is stable in both cases. Fig.~\ref{fig:example2} shows the largest inter-sample interval bound $h>0$ for which (\ref{eq:stabLMI}) can be verified with $X\geq 0$ and $Y\geq 0$ (blue-star 
$T_z=0.2$, green-plus $T_z=0.05$), and with $X\geq 0$ and $Y=0$ fixed (red-circle $T_z=0.2$, grey-diamond $T_z=0.05$). The latter corresponds to consideration of the gain-bound IQC only in the analysis. For this example, it is clear that accounting for the input feedforward passivity IQC from Lemma~\ref{iqc:pass} leads to much less conservative results. The reduction in conservativeness is more significant for the smaller value of $T_z$.   


\section{Conclusion} \label{sec:conc}

A model for a class of asynchronous sampled-and-hold operators is proposed and characterized by a family of IQCs. The model is parametrized by bounds on the uncertain inter-sample interval and the asynchrony between input sample events and zero-order-hold output update events. In principle,  the IQC representation can be used to devise robust stability and performance certificates for sampled-data networks of dynamical systems, as illustrated here for a feedback network with one asynchronous sample-and-hold link. For more complicated networks, it is expected that the IQC framework will facilitate exploitation of network structure as in~\cite{Kao09}, for example, and the consideration of non-linearities, such as link quantization as in~\cite{FabbroThesis,FabbroIFAC}, for example, and in the sub-system dynamics more generally. Within a networked systems context, additional IQC may be also used to encode known bounds on relationships between related links, as in the case of a common source and sampling sequence, but distinct hold update sequences, as considered in~\cite{CanFabKao18} for example. Ongoing work is continuing in these directions.


\addtolength{\textheight}{0cm}   








\end{document}